\documentclass[12pt]{amsart}
\usepackage{amscd}
\def\NZQ{\mathbb}               

\def\ZZ{{\NZQ Z}}
\def\RR{{\NZQ R}}

\def\frk{\mathfrak}               

\def\pp{{\frk p}}

\def\mm{{\frk m}}

\def\Phi{{\frk n}}
\def\Phi{{\frk N}}
\def\opn#1#2{\def#1{\operatorname{#2}}} 
\opn\chara{char} \opn\length{\ell} \opn\pd{pd} \opn\rk{rk}
\opn\projdim{proj\,dim} \opn\injdim{inj\,dim} \opn\rank{rank}
\opn\depth{depth} \opn\grade{grade} \opn\height{height}
\opn\embdim{emb\,dim} \opn\codim{codim} \opn\codepth{codepth}

\opn\Tr{Tr} \opn\bigrank{big\,rank}  \opn\Bor{Bor}
\opn\superheight{superheight}\opn\lcm{lcm}
\opn\trdeg{tr\,deg}%
\opn\reg{reg} \opn\lreg{lreg} \opn\ini{in}
\opn\div{div} \opn\Div{Div} \opn\cl{cl} \opn\Cl{Cl}
\opn\Spec{Spec} \opn\Supp{Supp} \opn\supp{supp} \opn\Sing{Sing}
\opn\Ass{Ass}
\opn\Ann{Ann} \opn\Rad{Rad} \opn\Soc{Soc}
\opn\Ker{Ker} \opn\Coker{Coker} \opn\Am{Am} \opn\Hom{Hom}
\opn\Tor{Tor} \opn\Ext{Ext} \opn\End{End} \opn\Aut{Aut}
\opn\id{id}

\opn\nat{nat}
\opn\pff{pf}
\opn\Pf{Pf} \opn\GL{GL} \opn\SL{SL} \opn\mod{mod} \opn\ord{ord}
\opn\aff{aff} \opn\con{conv} \opn\relint{relint} \opn\st{st}
\opn\lk{lk} \opn\cn{cn} \opn\core{core} \opn\vol{vol}
\opn\link{link} \opn\star{star} \opn\sort{sort}
\opn\gr{gr}

\def\pot#1#2{#1[\kern-0.28ex[#2]\kern-0.28ex]}

\opn\dirlim{\underrightarrow{\lim}}
\opn\inivlim{\underleftarrow{\lim}}
\let\union=\cup

\let\to=\rightarrow

\def\Implies{\ifmmode\Longrightarrow \else
     \unskip${}\Longrightarrow{}$\ignorespaces\fi}
\def\implies{\ifmmode\Rightarrow \else
     \unskip${}\Rightarrow{}$\ignorespaces\fi}
\def\iff{\ifmmode\Longleftrightarrow \else
     \unskip${}\Longleftrightarrow{}$\ignorespaces\fi}

\let\:=\colon
\newtheorem{Theorem}{Theorem}[section]
\newtheorem{Lemma}[Theorem]{Lemma}
\newtheorem{Corollary}[Theorem]{Corollary}
\newtheorem{Proposition}[Theorem]{Proposition}
\newtheorem{Remark}[Theorem]{Remark}

\newtheorem{Example}[Theorem]{Example}
\newtheorem{Examples}[Theorem]{Examples}

\let\epsilon\varepsilon
\let\phi=\varphi
\let\kappa=\varkappa
\def\qed{\ifhmode\textqed\fi
   \ifmmode\ifinner\quad\qedsymbol\else\dispqed\fi\fi}
\def\textqed{\unskip\nobreak\penalty50
    \hskip2em\hbox{}\nobreak\hfil\qedsymbol
    \parfillskip=0pt \finalhyphendemerits=0}
\def\dispqed{\rlap{\qquad\qedsymbol}}

\opn\dis{dis}
\def\pnt{{\raise0.5mm\hbox{\large\bf.}}}

\begin{document}

\title{Equidimensional and Unmixed Ideals of Veronese Type}
\author{ Marius Vladoiu }
\thanks{The author was supported by CNCSIS and the Ceres program 4-131/2004 of the Romanian
Ministery of Education and Research}

\address{Marius Vladoiu, Faculty of Mathematics, University of Bucharest, Str. Academiei 14,
Bucharest, 010014, Romania} \email{vladoiu@gta.math.unibuc.ro, marius\_alg@yahoo.com}

\begin{abstract}
This paper was motivated by a problem left by Herzog and Hibi, namely to classify all unmixed polymatroidal ideals. In the particular case of polymatroidal ideals corresponding to discrete polymatroids of Veronese type, i.e ideals of Veronese type, we give a complete description of the associated prime ideals and then, we show that such an ideal is unmixed if and only if it is Cohen-Macaulay. We also give for this type of ideals equivalent characterizations for being equidimensional. 
\end{abstract}

\maketitle

\section*{Introduction}

   In this paper, we give a partial answer to a problem left by Herzog and Hibi in \cite{HH1}. More precisely, they classify all Cohen-Macaulay polymatroidal ideals and, since a Cohen-Macaulay ideal is unmixed, they leave as a problem the classification of all unmixed polymatroidal ideals. Our main result gives a complete answer for this problem in the case of polymatroidal ideals corresponding to discrete polymatroids of Veronese type. A {\em polymatroidal ideal} is the ideal generated by the monomials corresponding to the set of bases of a discrete polymatroid. We recall that an ideal $J$ of a noetherian ring $R$ is said to be {\em unmixed} if the associated prime ideals of $R/J$ are the minimal prime ideals of $J$. In general, an ideal in a noetherian ring $R$ is called {\em equidimensional} if for all its minimal primes $\pp$, $\dim(R/\pp)$ is the same number.
   
   The paper is organized as follows. First of all, in Section $1$ we recall some basic facts about discrete polymatroids (see, e.g. \cite{HH}) and we focus on a particular class, namely the class of discrete polymatroids which satisfy the strong exchange property. Herzog et al.(\cite{HHV}) showed that the set of bases of a discrete polymatroid satisfying the strong exchange property is, up to an affinity, of Veronese type. The polymatroidal ideal corresponding to a discrete polymatroid of Veronese type we call it {\em ideal of Veronese type}, and denote it by $I_{d;a_1,a_2,\ldots,a_n}$, a subset of $S=K[x_1,\ldots,x_n]$, where $d,a_1,\ldots,a_n$ are nonnegative integers with $\sum_{i=1}^n a_i\geq d$. After renumbering the variables, we may assume that $d\geq a_1\geq a_2\geq\ldots\geq a_n$. With respect to this assumption, our first result in Section $2$ says that the radical of an ideal of Veronese type is squarefree strongly stable. Using this, we obtain the main theorem of this section, Theorem \ref{essen}, which states that an ideal of Veronese type $I$ is equidimensional if and only if $\sqrt{I}$ is Cohen-Macaulay. In particular, we obtain some combinatorial information about the minimal generators of $\sqrt{I}$, which allows us to give a new equivalent description of the  equidimensionality of $I$, cf. Corollary \ref{guso}. We show also that equidimensional polymatroidal ideals corresponding to discrete polymatroids which satisfy strong exchange property are those whose radical is a principal ideal. In the last section, we give in Proposition \ref{radoi} a complete description of the associated prime ideals of $S/I$, where $I$ is an ideal of Veronese type. Finally, we prove the main result of the paper, Theorem \ref{vlad}, which asserts that an ideal of Veronese type, $I$, is unmixed if and only if $I$ is unmixed and equidimensional if and only if $I$ is Cohen-Macaulay. In particular, this shows that unmixed implies equidimensional for this type of polymatroidal ideals and, it follows from \cite[Theorem 3.2]{HH1}, that the unmixed ideals of Veronese type are precisely the principal ideals, the Veronese ideals, and the squarefree Veronese ideals.
 
 I am very grateful to Professor J\" urgen Herzog for many helpful discussions and comments during my stay at the University Duisburg-Essen. Also, I wish to thank for the warm hospitality of School of Mathematical Sciences from Lahore, Pakistan, where I could finish this paper.

\section{Review on Discrete Polymatroids}

\bigskip

  Fix an integer $n > 0$ and set $[n] = \{ 1, 2, \ldots, n\}$. The canonical basis vectors of $\RR^n$ will be denoted by $\varepsilon_1, \ldots, \varepsilon_n$. Let $\ZZ_+^n$ denote the set of those vectors $u =(u_1, \ldots, u_n) \in
\ZZ^n$ with each $u_i \geq 0$. For a vector $u =(u_1, \ldots, u_n) \in \ZZ_+^n$, the {\em modulus} of $u$ is the number $|u| =  \sum_{i=1}^n u_i$. If $u=(u_1,\ldots,u_n)$ and $v=(v_1,\ldots,v_n)$ are two vectors belonging to $\ZZ_+^n$, then we write $u\leq v$ if all components $v_i-u_i$ of $v-u$ are nonnegative, and moreover, write $u<v$ if $u\leq v$ and $u\neq v$. 

  A {\em discrete polymatroid} (see \cite{HH}) on the ground set $[n]$ is a
non-empty finite subset $P\subset \ZZ^n_+$ satisfying

\begin{enumerate}
\item[(1)] if $u  \in P$ and $v  \in \ZZ_+^n$ with $v \leq u$, then $v \in P$; 
\item[(2)] if $u = (u_1, \ldots, u_n) \in P$ and $v =(v_1, \ldots, v_n) \in P$ with $|u| < |v|$, then there is $i \in
[n]$ with $u_i < v_i$ such that $u + \varepsilon_i \in P$.
\end{enumerate}
A {\em base} of $P$ is a vector $u \in P$ such that $u < v$ for no $v \in P$. We denote the set of bases by $B(P)$. It is an easy consequence of the definition of a discrete polymatroid that all elements in $B(P)$ have the same modulus. This common number is called the {\em rank} of $P$ and denote it by $\rank(P)$.

According to  \cite[Theorem 2.3]{HH} the set of bases of a discrete polymatroid can be characterized by the {\em exchange property}: a subset $B\subset \ZZ_+^n$ of vectors of the same modulus, is the set of bases of a discrete polymatroid if and only if  for all $u, v \in B$ such that  $u_i > v_i$ for some $i$,  there exists  $j \in [n]$ with $u_j < v_j$ such that $u - \varepsilon_i + \varepsilon_j \in B$.

We say that a discrete polymatroid satisfies the {\em strong exchange property}, if for all $u,v\in  B(P)$ and for any $i,j\in [n]$ such that $u_i>v_i$ and $u_j<v_j$, one has that $u -\varepsilon_i + \varepsilon_j \in B(P)$.

An important characterization for the discrete polymatroids is given by the rank function.  Let $P$ be a discrete polymatroid. The {\em ground set rank function} of $P$ is the  function $\rho_P: 2^{[n]} \to \ZZ_{+}$ defined by setting
\[
\rho_P(A) = \max\{ v(A) \: v \in P \}
\]
for all $\emptyset \neq A \subset [n]$ together with $\rho_P(\emptyset) = 0$.

This function is {\em nondecreasing}, i.e., if $A \subset B \subset [n]$, then $\rho_P(A) \leq \rho_P(B)$, and is {\em
submodular}, i.e.,
\[
\rho_P(A) + \rho_P(B) \geq \rho_P(A \cup B) + \rho_P(A \cap B),
\]
for all $A, B \subset [n]$.


In the following examples we give a class of discrete polymatroids which satisfy the strong exchange property.

\begin{Examples}
\label{examples1}{\em (a) Let $a_1,\ldots,a_n$ and $d$ be nonnegative integers such that $a_i\leq d$ for $i=1,\ldots,n$. The non-empty subset $B$ of $\ZZ^n_+$ (non-empty means that $\sum_{i=1}^n a_i\geq d$) given by
\[
B:=\{u\in \ZZ^n_+: |u|=d, \quad 0\leq u_i\leq a_i \}
\]
is the set of bases of a discrete polymatroid of rank $d$, called {\em discrete polymatroid of Veronese type} and denoted by $P_{d;a_1,\ldots,a_n}$. It satisfies the strong exchange property and its ground set rank function $\rho$ is given by:
\begin{eqnarray}
\label{rank}
 \rho(A)=\min(\sum_{i\in A} a_i, d),  \text{ for any }
A\subset [n].
\end{eqnarray}
 In particular $\rho([n])=d$.

(b) The set $\{(2,1,1),(1,2,1),(1,1,2)\}$ is the set of bases of a discrete polymatroid of rank $4$ which satisfies the strong exchange property, but it is not of Veronese type. Its ground set rank function is given by $\rho(\{1\})=\rho(\{2\})=\rho(\{3\})=2$,
$\rho(\{1,2\})=\rho(\{1,3\})=\rho(\{2,3\})=3$ and $\rho([3])=4$.

(c) A first example of a discrete polymatroid which does not satisfy the strong exchange property arises in dimension $4$. The set $\{(1,0,1,0),(0,1,1,0),(0,1,0,1),\\ (1,0,0,1)\}$ is the set of bases of a discrete polymatroid of rank $2$, but does not satisfy the strong exchange property. Its ground set rank function is given by $\rho(\{1\})=\rho(\{2\})=\rho(\{3\})=\rho(\{4\})=1$,
$\rho(\{1,2\})=\rho(\{3,4\})=1$,
$\rho(\{1,3\})=\rho(\{1,4\})=\rho(\{2,3\})=\rho(\{2,4\})=2$,
$\rho(\{1,2,3\})=\rho(\{1,2,4\})=\rho(\{1,3,4\})=\rho(\{2,3,4\})=2$,
$\rho([4])=2$.}
\end{Examples}

The following theorem (see \cite[Theorem 1.1]{HHV}) , which will be needed later, shows the relationship between discrete polymatroids which satisfy strong exchange property and discrete polymatroids of Veronese type:

\begin{Theorem}
\label{raul} Let $P$ be discrete polymatroid which satisfies the strong exchange property. Then $B(P)$ is isomorphic to the set of bases of a polymatroid of Veronese type.
\end{Theorem}

We fix the following notation: let $K$ be a field and $S=K[x_1,\ldots,x_n]$ the polynomial ring in $n$ variables over $K$ with each $\deg x_i=1$. For $I\subset S$ a monomial ideal we denote by $G(I)$ its unique set of minimal generators. For a monomial $u\in S$ we define by $\max(u)=\max\{i| \ x_i \text{ divides } u\}$ and by $\min(u)=\min\{i| \ x_i \text{ divides } u\}$. For a subset $A=\{i_1,\ldots,i_k\}$ of $[n]$ we denote by $x_A$ the monomial $x_{i_1}\cdots x_{i_k}$, and by $P_A$ the prime ideal of $S$ generated by the variables whose index is in $A$, i.e. the ideal
$(x_{i_1},\ldots,x_{i_k})$. 

We recall that a monomial ideal generated by monomials corresponding to the set of bases of a discrete polymatroid $P$ is called a {\em polymatroidal ideal}, and denoted $I(P)$. In other words, $I(P)$ is generated by all monomials $x^v$ with $v\in B(P)$.

In order to shorten the name, we shall use instead of
polymatroidal ideal corresponding to a discrete polymatroid with
strong exchange property, the name {\em strong polymatroidal
ideal}, and instead of polymatroidal ideal corresponding to a
discrete polymatroid of Veronese type, the name {\em ideal of
Veronese type}.

\begin{Remark}
\label{figo}{\em For any polymatroidal ideal $I=I(P)\subset S$, the
radical $\sqrt{I}$ of $I$ is minimally generated by all monomials
$x_A$, where  $A\subset [n]$ satisfies $\rho_P(A)=\rank(P)$ and is minimal
with respect to inclusion. }
\end{Remark}

\section{Classification of equidimensional ideals of Veronese type}

Let $I\subset S$ be a monomial ideal. A {\em vertex cover} of $I$
is a subset $W$ of $\{x_1,\ldots,x_n\}$ such that each $u\in G(I)$
is divided by some $x_i\in W$. Such a vertex cover $W$ is called
{\em minimal} if no proper subset of $W$ is a vertex cover of $I$.
The ideals generated by the minimal vertex covers of $I$ are
exactly the minimal prime ideals of $I$. We recall, that in
general, an ideal in a Noetherian ring $S$ is called {\em
equidimensional} if for all its minimal primes $\pp$, $\dim S/\pp$
is the same number. It is easy to see that a monomial ideal $I$ is
equidimensional if all minimal vertex covers of $I$ have the same
cardinality.

Throughout this section we consider $I$ an ideal of Veronese type,
i.e. there exists nonnegative numbers $a_1,\ldots,a_n$ and $d$
with $a_i\leq d$ for $i=1,\ldots,n$ such that
$I=I(P_{d;a_1,\ldots,a_n})$, and denote it for the rest of paper by
$I_{d;a_1,\ldots,a_n}$. We may assume, after renumbering the variables, that $a_1\geq a_2\geq \ldots\geq a_n$. 

It follows from (\ref{rank}) that $G(\sqrt{I})$ is the set of all
monomials
$$x_{i_1}\cdots x_{i_k}\quad\text{with}\quad k\geq 1,\quad
1\leq i_1<\ldots< i_k\leq n,
$$  such that
\begin{eqnarray}
\label{rad} a_{i_1}+\cdots +a_{i_k}\geq d \quad\text{and}\quad
a_{i_1}+\cdots +a_{i_{k-1}}<d.
\end{eqnarray}

Obviously if $I$ is a monomial ideal, $I$ is equidimensional if
and only if $\sqrt{I}$ is equidimensional.

We may also assume that $d> a_1\geq a_2\geq\ldots\geq a_n\geq 1$, so $I$ is an ideal of $S$. Indeed, if $a_{i+1}=\ldots=a_n=0$ for some $1<i<n$ and $a_i> 0$ then we take $I=I_{d;a_1,\ldots,a_i}$ and consider it as an ideal of $K[x_1,\ldots,x_i]$. If $a_1=\ldots=a_n=d$, then $I$ is a Veronese ideal, so it is Cohen-Macaulay and hence equidimensional. Otherwise, if $a_1=\ldots=a_l=d$ for some $1\leq l<n$, then $x_1,\ldots,x_l\in G(\sqrt{I})$ and hence any minimal vertex cover of $G(\sqrt{I})$ contains $\{x_1,\ldots,x_l\}$, so $\sqrt{I}$ is equidimensional if and only if the Veronese type ideal $I_{d;a_{l+1},\ldots,a_n}\subset K[x_{l+1},\ldots,x_n]$ is equidimensional. But now $d>a_{l+1}\geq\ldots \geq a_n$, and therefore we can consider for the rest of this section that $I=I_{d;a_1,\ldots,a_n}$, with $d>a_1\geq\ldots\geq a_n\geq 1$.

The following lemma is easy but crucial for this section.
\begin{Lemma}
\label{ionescu} Let $I$ be an ideal of Veronese type. Then, 
$\sqrt{I}$ is squarefree strongly stable. 
\end{Lemma}

\begin{proof}
Let $u\in G(\sqrt{I})$. Then, according to (\ref{rad})
$u=x_{i_1}\cdots x_{i_k}$ with $1\leq i_1<\ldots <i_k\leq n$, such
that $a_{i_1}+\cdots +a_{i_k}\geq d$ and $a_{i_1}+\cdots
+a_{i_{k-1}}<d$. Let $j\in\{i_1,\ldots,i_k\}$ and $l<j$ with
$l\notin\{i_1,\ldots,i_k\}$. Since $a_l\geq a_j$ we have that
$a_{i_1}+\cdots +a_{i_k}-a_j+a_l\geq a_{i_1}+\cdots +a_{i_k}\geq d$,
hence $u\cdot\frac{x_l}{x_j}\in\sqrt{I}$, so $\sqrt{I}$ is strongly
stable.
\end{proof}

\begin{Remark}
{\em In general one can associate to any discrete polymatroid $P\subset\ZZ_+^n$, with the rank function $\rho$, a simplicial complex (see for example \cite{BH},\cite{V}) $\Delta$ on the ground set $[n]$, as follows: $F\subset [n]$ is a face of $\Delta$ if $\rho(F)<\rho([n])$. In the particular case of discrete polymatroid of Veronese type of rank $d$, 
$\Delta$ is the collection of subsets $G\subset [n]$ such that $\sum_{i\in G}a_i<d$, hence it contains all the vertices (see the definition \cite[Chapter 5]{BH}), and the Stanley-Reisner ideal of $\Delta$ is $I_{\Delta}=\langle x_F| F\notin\Delta\rangle=\sqrt{I}$.

To any simplicial complex $\Delta$ one can associate its Alexander
dual $\Delta^{\vee}$, which is the simplicial complex with
Stanley-Reisner ideal $I_{\Delta^{\vee}}=\langle x_{G^c}|
G\in\Delta\rangle$ (see for example \cite[Lemma 2.1]{He}). Here by
$G^c$ we denote the set $[n]\setminus G$. In general for any
squarefree monomial ideal $J\subset S$ one can define its
Alexander dual $J^{\vee}$. Because $J$ is squarefree, then there
exists a simplicial complex $\Gamma$ such that $J=I_{\Gamma}$.
Then, we set $J^{\vee}:=I_{\Gamma^{\vee}}$. }
\end{Remark}

We recall

\begin{Lemma}
\label{zidane}
 Let $J\subset S$ be a squarefree strongly stable ideal. Then $J^{\vee}$
 is again a squarefree strongly stable ideal.
\end{Lemma}

\begin{proof}

As in the previous remark, let $\Gamma$ be the simplicial complex
such that $J=I_{\Gamma}$, i.e. $J=\langle x_F|
F\not\in\Gamma\rangle$. Let us first observe that
$G(I_{\Gamma^{\vee}})=\{x_{G^c}| G \text { is a facet of }
\Gamma\}$. Consider now $m\in G(I_{\Gamma^{\vee}})$, i.e.
$m=x_{j_1}\cdots x_{j_{n-k}}$, where
$\{j_1,\ldots,j_{n-k}\}=[n]\setminus\{i_1,\ldots,i_k\}$ with
$\{i_1,\ldots,i_k\}$ a facet of $\Gamma$. We have to show that for
any $j\in \{j_1,\ldots,j_{n-k}\}$ and for any $i<j$ such that
$i\not\in\{j_1,\ldots,j_{n-k}\}$ the monomial
$m\cdot\frac{x_i}{x_j}\in J^{\vee}$. Since
$m\cdot\frac{x_i}{x_j}=x_{D^c}$, where
$D=(\{i_1,\ldots,i_k\}\setminus\{i\})\cup\{j\}$, if we prove that
$D\in\Gamma$, then we are done. Assume on the contrary that
$D\not\in\Gamma$. Then $u=x_D=x_{i_1}\cdots x_{i_k}(x_j/x_i)\in
J$, and because $J$ is strongly stable and $x_j|u$, $i<j$ and
$i\not\in\Supp(u)$ implies $u\cdot\frac{x_i}{x_j}\in J$, hence
$x_{i_1}\cdots x_{i_k}\in J$, or equivalently
$\{i_1,\ldots,i_k\}\not\in\Gamma$, a contradiction. Therefore,
$J^{\vee}$ is squarefree strongly stable.
\end{proof}

It follows then from Lemma \ref{ionescu} and Lemma \ref{zidane} that
for an ideal of Veronese type $I$, the ideals  $\sqrt{I}$ and
$(\sqrt{I})^{\vee}$ are squarefree strongly stable ideals. We
introduce now for a squarefree monomial ideal $J\in S$ the
following invariants: $$m(J)=\max\{\max(u): u\in G(J)\}$$ and
$$b(J)=\max\{\deg(u): u\in G(J)\}.$$
\begin{Lemma}
\label{samuel} Let $I=I_{d;a_1,\ldots,a_n}\subset S$ be an ideal of
Veronese type. Then
\begin{enumerate}
\item[(i)] $m(\sqrt{I})$ equals the maximum of the integers $i_k$
such that there exists $$1\leq i_1<\ldots<i_k\leq n \text{ with }
a_{i_1}+\ldots+a_{i_k}\geq d  \text{ and }
a_{i_1}+\ldots+a_{i_{k-1}}<d;$$ \item[(ii)] $b(\sqrt{I})=\max\{l|
\text{ such that there exists } 1\leq i_1<\ldots<i_l\leq n \text{
with } a_{i_1}+\ldots+a_{i_l}\geq d \text{ and }
a_{i_1}+\ldots+a_{i_{l-1}}<d\}$.
\end{enumerate}
\end{Lemma}

\begin{proof}
The assertions  follow immediately from the inequalities
(\ref{rad}).
\end{proof}

\begin{Remark}
{\em (a) For $i_k$ which equals $m(\sqrt{I})$ we can have different $k$'s. For
example if $I=I_{7;4,3,2,1,1}$, then
\[
G(\sqrt{I})=\{x_1x_2,x_1x_3x_4,x_1x_3x_5,x_2x_3x_4x_5\},
\]
so $m(\sqrt{I})=5$. We have two different sequences, whose maximal element
equals $5$, $1<3<5$ ( here $5=i_3$, so $k=3$) and $2<3<4<5$ (here
$5=i_4$, so $k=4$).

(b) One could expect that $b(\sqrt{I})=\max\{k| \text{ such that there exists
} 1\leq i_1<\ldots<i_k=m \text{ with } a_{i_1}+\ldots+a_{i_k}\geq d
\text{ and } a_{i_1}+\ldots+a_{i_{k-1}}<d\}$, which is true in the
previous example. This is not true in general, as we can see from the
following. Let $I=I_{11;7,4,3,2,2,1}$, then
\[
G(\sqrt{I})=\{x_1x_2,x_1x_3x_4,x_1x_3x_5,x_1x_3x_6,x_1x_4x_5,x_2x_3x_4x_5\}.
\]
We observe that $m(\sqrt{I})=6$, and if our assumption would be true, then
$b(\sqrt{I})$ should be $3$, a contradiction since $b(\sqrt{I})=4$.

(c) For a strongly stable ideal $J\subset S$, let us denote by $\Bor(J)$ the minimal set of Borel generators of $J$. For example, for the ideals above $$\Bor(\sqrt{I_{7;4,3,2,1,1}})=\{x_1x_2,x_1x_3x_5,x_2x_3x_4x_5\}$$ and $$\Bor(\sqrt{I_{11;7,4,3,2,2,1}})=\{x_1x_2,x_1x_3x_6,x_1x_4x_5,x_2x_3x_4x_5\}.$$ Then, it is easy to see that for a squarefree strongly stable ideal $J$, $m(J)=\max\{\max(u): u\in \Bor(J)\}$ and $b(J)=\max\{\deg(u): u\in \Bor(J)\}$. 
}
\end{Remark}

\begin{Theorem}
\label{essen} Let $I$ be an ideal of Veronese type. Let
$m=m(\sqrt{I})$ and $b=b(\sqrt{I})$. Then the following are
equivalent:
\begin{enumerate}
\item[(a)] $I$ is equidimensional, \item[(b)] $\sqrt{I}$ is
Cohen-Macaulay, \item[(c)] There is a unique Borel generator of
degree $b$ of $\sqrt{I}$, namely
$$x_{m-b+1}\cdots x_{m},$$ and for any $u\in
 G(\sqrt{I})$ we have $\max(u)-\deg(u)\leq m-b$.
\end{enumerate}
Moreover, if $I$ is equidimensional, then the cardinality of any
minimal vertex cover of $I$ is $m-b+1$.
\end{Theorem}

\begin{proof}
 We shall prove that $(a)\Leftrightarrow (b)$ and $(b)\Leftrightarrow
 (c)$.

$(a)\Leftrightarrow (b)$: If the ideal $\sqrt{I}$ is
Cohen-Macaulay, then $\sqrt{I}$ is equidimensional, so $I$ is
equidimensional. For the converse assume that $I$ is
equidimensional. Then $\sqrt{I}=I_{\Delta}$ is equidimensional.
According to the \cite[Lemma 2.2]{He}, this means that
$I_{\Delta^{\vee}}$ is generated in one degree. By Lemma
\ref{zidane}, $I_{\Delta^{\vee}}$ is squarefree strongly stable,
hence it has linear quotients. It follows from \cite[Lemma
4.1]{CH} that $I_{\Delta^{\vee}}$ has a linear resolution over
$S$. Now, applying Eagon-Reiner theorem \cite{ER}, we obtain that
$\sqrt{I}=I_{\Delta}$ is
 a Cohen-Macaulay ideal.

(b)\implies (c):  Herzog and Srinivasan compute in
\cite[Proposition 4.1]{HS} the codepth and codimension of $S/J$,
with $J$ a squarefree strongly stable ideal:
$$\codim S/J=\max\{\min(u): u\in G(J)\},$$ $$\codepth
S/J=\max\{\max(u)-\deg(u):u\in G(J)\}+1.$$ It is easy to see that
in these formulas $G(J)$ may be replaced by $\Bor(J)$.
Therefore (b) implies  that
\begin{eqnarray}
\label{equal} \quad \max\{\min(u):u\in
\Bor(\sqrt{I})\}=\max\{\max(u)-\deg(u):u\in \Bor(\sqrt{I})\}+1.
\end{eqnarray}

Let  $\Bor(\sqrt{I})=\{u_1,\ldots,u_r\}$, and let
$i_0,j_0\in\{1,\ldots,r\}$ such that
\begin{eqnarray}
\label{depth} \max(u)-\deg(u)\leq \max(u_{j_0})-\deg(u_{j_0})
\quad \text{for all}\quad u\in \Bor(\sqrt{I}),
\end{eqnarray}
 and
\begin{eqnarray}
\label{dim} \min(u)\leq \min(u_{i_0}) \quad \quad \text{for
all}\quad u\in \Bor(\sqrt{I}).
\end{eqnarray}
By (\ref{equal}) we get
\[
\min(u_{i_0})\leq \max(u_{i_0})-\deg(u_{i_0})+1\leq
\max(u_{j_0})-\deg(u_{j_0})+1=\min(u_{i_0}),
\]
and hence  $\min(u_{i_0})=\max(u_{i_0})-\deg(u_{i_0})+1$. We
conclude that $u_{i_0}=x_t\cdots x_{t+b'-1}$, where we denoted by
$t=\min(u_{i_0})$ and $b'=\deg(u_{i_0})$.

It remains to be shown that $b'=b$ and $t=m-b+1$. Suppose there is
an element  $u\in \Bor(\sqrt{I})$ such that $\deg(u)\geq b'$.

Assume first that $\deg(u)=b'$, i.e. $u=x_{l_1}\cdots x_{l_{b'}}$
for some $l_1<\ldots <l_{b'}$ in $[n]$. (\ref{dim}) and
(\ref{depth}) implies that $l_1\leq t$ and $l_{b'}\leq t+b'-1$.
Therefore $l_s\leq t+s-1$, for any $1\leq s\leq b'$, so $u$
belongs to the strongly stable ideal generated by $u_{i_0}$, a
contradiction with $u\in \Bor(\sqrt{I})$.

Let now $k=\deg(u)>b'$ and $u=x_{l_1}\cdots x_{l_k}$ for a
sequence $l_1<\ldots <l_k$ in $[n]$. (\ref{dim}) implies that
$l_1\leq t$ and from (\ref{depth}) it follows that $l_k\leq
t+k-1$, therefore $l_{b'}\leq l_k-(k-b')\leq t+b'-1$. It is an
easy consequence of previous inequality and to our non-descending
sequence $a_1\geq\ldots\geq a_n$ that
$a_{l_1}+\ldots+a_{l_{b'}}\geq a_t+\ldots+a_{t+b'-1}\geq d$ (the
last inequality is due to $u_{i_0}\in \Bor(\sqrt{I})\subset
G(\sqrt{I})$). Since $u\in \Bor(\sqrt{I})$ we have that
$d>a_{l_1}+\ldots+a_{l_{k-1}}\geq a_{l_1}+\ldots+a_{l_{b'}}\geq
d$, a contradiction.

Hence we proved our claim is true. Therefore $b'=b$ and
$\deg(u)<b$, $\forall u\in \Bor(\sqrt{I})\setminus\{u_{i_0}\}$. We
have also $\max(u)-\deg(u)\leq t-1$, $\min(u)<t$, $\forall u\in
\Bor(\sqrt{I})\setminus\{u_{i_0}\}$ (otherwise, if $\min(u)=t$ we
obtain $u=u_{i_0}$). It follows that $\max(u)\leq
t+\deg(u)-1<t+b-1$, $\forall u\in
\Bor(\sqrt{I})\setminus\{u_{i_0}\}$. Hence $m=t+b-1$,
$u_{i_0}=x_{m-b+1}\cdots x_m$ is the unique Borel generator of
degree $b$, and the rest of $(c)$ also holds.

$(c)\Rightarrow (b)$: It is straightforward to check that
condition (c) implies $\codepth S/\sqrt{I}=\codim S/\sqrt{I}$.
Hence $\sqrt{I}$ is a Cohen-Macaulay ideal, and so $I$ is
equidimensional.

If the equivalent conditions hold, then all the facets of $\Delta$
have the same cardinality. All the subsets of cardinal $b$ of
$[m]$ are non-faces of $\Delta$ since $x_{m-b+1}\cdots x_m$ is a
Borel generator of $I_{\Delta}$. Hence the cardinality of a facet
of $\Delta$ is $\leq b-1$. Since $x_{m-b+1}\cdots x_{m-1}\not\in
I_\Delta$ it follows that $\{m-b+1,\ldots,m-1\}$ is a face of
$\Delta$ of cardinality $b-1$, and hence a facet. Therefore, all
the facets have cardinality $b-1$. Consequently, all minimal vertex covers have cardinality $m-b+1$.
\end{proof}

\begin{Remark}
\label{cesar}{\em Both conditions in statement (c) of Theorem
\ref{essen}  are needed.

(a) If $I_1=I_{9;7,3,3,2,1}$, then
$$G(\sqrt{I_1})=\{x_1x_2,x_1x_3,x_1x_4,x_2x_3x_4x_5\}$$ and
$$\Bor(\sqrt{I_1})=\{x_1x_4,x_2x_3x_4x_5\}.$$ Since $m=5$ and $b=4$,
we see that $u=x_2x_3x_4x_5$ is the unique Borel generator of
maximal degree with $\max(u)=5$.  On the other hand,
$\max(x_1x_4)-\deg(x_1x_4)=4-2=2>1=5-4=m-b$. Therefore, according
to the theorem, $I_1$ is not equidimensional.  Indeed,
$\{x_1,x_2\}$ and $\{x_2,x_3,x_4\}$ are minimal vertex covers of
$\sqrt{I_1}$ of different cardinality.

(b) If $I_2=I_{8;5,5,4,3,1,1}$, then
$$G(\sqrt{I_2})=\{x_1x_2,x_1x_3,x_1x_4,x_2x_3,x_2x_4,x_3x_4x_5,x_3x_4x_6\}$$
and $$\Bor(\sqrt{I_2})=\{x_2x_4,x_3x_4x_6\}.$$ Furthermore, we
have $m=6$, $b=3$, $\max(x_2x_4)-\deg(x_2x_4)=2<3=6-3$ and
$\max(x_3x_4x_6)-\deg(x_3x_4x_6)=3\leq 3=6-3$. So for any
$u\in\Bor(\sqrt{I_2})$ we have $\max(u)-\deg(u)\leq m-b$. But
$I_2$ is not equidimensional, according to the theorem, since the
Borel generator in maximum degree is $x_3x_4x_6$ instead of
$x_4x_5x_6$. Indeed, one can check that $\{x_1,x_2,x_3\}$ and
$\{x_1,x_2,x_5,x_6\}$ are minimal vertex covers of $\sqrt{I_2}$ of
different cardinality. }
\end{Remark}

\begin{Corollary}
\label{guso}
Let $I=I_{d;a_1,\ldots, a_n}$ be an ideal of Veronese type. Then, the ideal $I$ is equidimensional if and only if there exists a pair $(p,l)$ with $p\geq l$, which is maximal with respect to the partial order given by componentwise comparison, and such that
\begin{enumerate}
\item[(i)] $a_{p-l+1}+\ldots+a_{p}\geq d$, and $a_{p-l+1}+\ldots+a_{p-1}<d$, and 
\item[(ii)] for any sequence $i_1<\ldots<i_k$ in $[n]$ with $a_{i_1}+\ldots+a_{i_k}\geq d$ and $a_{i_1}+\ldots+a_{i_{k-1}}<d$ we have $i_k-k\leq p-l$.
\end{enumerate}
\end{Corollary}

\begin{proof}
If $I$ is equidimensional, then it follows easily from Theorem \ref{essen} that the pair $(m,b)$ satisfies all the conditions from the corollary. 

For the converse, first we will show that any two pairs
$(p_1,l_1)$ and $(p_2,l_2)$ with the properties (i) and (ii) can be compared,
hence the name maximal is correct. Assume the contrary, then the only possibility (after a renumbering) for the two pairs to be incomparable with respect to our partial order is to have the following inequalities: $p_1<p_2$ and $l_1>l_2$. Then, $$d>a_{p_1-l_1+1}+\ldots+a_{p_1-1}\geq a_{p_2-l_2+1}+\ldots+a_{p_2}\geq d,$$ 
a contradiction (for the second inequality we use that $l_1-1\geq l_2$, $p_1<p_2$ and the non-increasing sequence $a_1\geq a_2\geq\cdots\geq a_n$). Therefore, the existence of a pair leads to existence of a maximal one.

For finishing the proof it is enough to show that the maximal pair $(p,l)$ from the hypotheses is equal to $(m,b)$. We have from the definition of $m$ and $b$ the following inequalities: $p\leq m$ and $l\leq b$. Arguing by contradiction,  suppose that $(p,l)<(m,b)$. This implies $p<m$ and $l\leq b$. Indeed, if $p=m$ and $l<b$, then from the definition of $b$ and $p$ we have the following inequalities: $$a_{m-l+1}+\ldots+a_{m}=a_{p-l+1}+\ldots+a_{p}\geq d$$ and
 $$a_{m-b+1}+\ldots+a_{m-1}<d,$$ which should be true simultaneously, a contradiction since $l<b$ and $a_1\geq a_2\geq\ldots\geq a_n$.

Hence $p<m$ and $l\leq b$. We have two cases to analyze: $l<b$ and
$l=b$.

First, let us suppose that $l<b$. From the definition of $b$, there exists a sequence $j_1<\cdots\ <j_b$ in $[n]$ such that $a_{j_1}+\ldots+a_{j_b}\geq d$, $a_{j_1}+\ldots+a_{j_{b-1}}<d$ and $j_b-b\leq p-l$. Since $(p,l)$ is maximal then $j_1<j_b-b+1\leq p-l+1$. In order to have both $a_{j_1}+\ldots+a_{j_{b-1}}<d$ and $a_{p-l+1}+\ldots+a_{p}\geq d$ we must have $j_l\geq p+1$. It follows then, that $j_b\geq j_l+b-l\geq p+1+b-l$, hence $j_b-b\geq p-l+1$, a contradiction since $j_b-b\leq p-l$.

Assume now that $l=b$. From the definition of $m$, there exists a sequence $i_1<\cdots <i_k=m$ in $[n]$ such that $a_{i_1}+\ldots+a_{i_k}\geq d$, $a_{i_1}+\ldots+a_{i_{k-1}}<d$. In addition, according to the hypotheses on the pair $(p,l)$, $m-k\leq p-l$. Since $m>p$ and $b=l\geq k$ we obtain that $m-k>p-l$, a contradiction. Hence, our discussion shows that in both cases we get a contradiction, so $(p,l)=(m,b)$. Applying now Theorem \ref{essen}, we obtain that $I$ is equidimensional. 
\end{proof}

\begin{Remark}
{\em (a) It follows from the proof of the corollary that the maximality of the pair $(p,l)$ with respect to the partial order given by the componentwise comparison is induced only by condition (i). 

(b) Corollary \ref{guso} provides an useful tool to show that an ideal of Veronese type is {\em not} equidimensional. The strategy is to compute first the maximal pair $(p,l)$ with the property (i) and then find a generator which does not satisfy property (ii). For example, consider the ideal of Veronese type $I=I_{15;9,6,4,3,2,2,1,1}$. We get succesively: $9+6=15$ implies that $(p_0,l_0)=(2,2)$, then $6+4<15$, $6+4+3<15$, $6+4+3+2=15$ implies that $(p_1,l_1)=(5,4)$, $4+3+2+2<15$, $4+3+2+2+1<15$, $4+3+2+2+1+1<15$. Therefore, the maximal pair $(p,l)$ is $(5,4)$ and $p-l=1$. It is easy to see that $x_1x_3x_6\in G(\sqrt{I})$ and, since $6-3=3>1=p-l$, the second condition (ii) of the Corollary \ref{guso} is not fulfilled, so we can conclude that $I$ is not equidimensional.}
\end{Remark}

Let $P$ be a discrete polymatroid which satisfies the strong
exchange property, with the rank function $\rho$. Then according to
the Theorem \ref{raul}, $B(P)$ is isomorphic to the set of bases of
a polymatroid of Veronese type. In fact, the isomorphism used in the
proof of the theorem is given by the translation
$\tau:\RR^n\longrightarrow\RR^n$,
$$\tau(v)=v-u_0, \quad \quad \forall v\in\RR^n,$$ where $u_0$ is
the integer vector whose $i$-th coordinate is equal to
$\rho([n])-\rho([n]\setminus\{i\})$. $\tau(B(P))$ is the set of
bases of a discrete polymatroid of Veronese type. If $u_0=0$, then
$P$ is of Veronese type. Eventually, after a renumbering of
variables, we may assume that
$\tau(B(P))=B(P_{d;a_1,\ldots,a_n})$, where $d\geq
a_1\geq\ldots\geq a_n$.

\begin{Corollary}
\label{owen} Let $I$ be a strong polymatroidal ideal which is not
of Veronese type. Then $I$ is equidimensional if and only if
$\sqrt{I}$ is principal.
\end{Corollary}

\begin{proof}
By the considerations preceding this corollary it follows that
$I=u_0J$, where $u_0\neq 1$ is a monomial and $J$ is a monomial
ideal. For any $x_k$ dividing $u$, $(x_k)$ is a minimal prime
ideal of $I$. Therefore $I$ is equidimensional if and only if all
minimal prime ideals of $I$, and hence of $\sqrt{I}$, have height
$1$. Let $P_i=(x_{j_i})\subset S$, for $1\leq i\leq r$, be these
minimal prime ideals. Since $\sqrt{I}=\bigcap_{i=1}^r P_i$ it
follows that $\sqrt{I}=(\prod_{i=1}^r x_{j_i})$. Hence $\sqrt{I}$
is principal. The converse implication  is obvious.
\end{proof}

\medskip

\section{Classification of unmixed ideals of Veronese type}

\bigskip

An ideal $I$ of a Noetherian ring $R$ is said to be {\em unmixed} if the associated prime ideals of $R/I$ are the minimal prime ideals of $I$. Some authors, for example the book of  Matsumura \cite{M}, require in addition that all minimal prime ideals have the same height. In the case of a polynomial ring this definition of unmixedness is equivalent to say that $I$ is equidimensional and has no embedded prime ideals.

The conditions ``unmixed'' and ``equidimensional'' are unrelated. For example, the ideal $(xy,xz)$ of $R=K[x,y,z]$ is unmixed but not equidimensional, while the ideal $(x^2,xy)$ is equidimensional but not unmixed.

Deciding whether an ideal is unmixed is in general a hard problem since one has to know all its associated prime ideals. In the particular case of ideals of Veronese type we can give a complete description of the associated prime ideals. As before, we may assume that our non-zero ideal of Veronese type is of the form $I=I_{d;a_1,\ldots,a_n}$ with $d\geq a_1\geq\ldots\geq a_n$. In case  $I$ is a prime ideal we have  $\Ass(S/I)=I$. In our case this happens only if $I=(x_1,\ldots,x_i)$. In our further discussions we therefore assume that  $d>1$. We may also assume that $a_n\geq 1$, because if $i$ is the least integer with $a_i=0$, then $I=I_{d;a_1,\ldots,a_{i-1}}S$ and $I$ is unmixed if and only if $I_{d;a_1,\ldots,a_{i-1}}\subset K[x_1,\ldots,x_{i-1}]$ is unmixed.

\begin{Proposition}
\label{radoi} Let $I=I_{d;a_1,\ldots,a_{n}}\subset S$ be an ideal of Veronese type with $d>1$ and $a_n\geq 1$, and $A$ a subset of $[n]$. Then 
$$P_A\in\Ass(S/I)\iff \sum_{i=1}^n a_i\geq d-1+|A|\quad \text{and}\quad \sum_{i\notin A} a_i\leq d-1.$$ 
 Moreover, if for all $i\in A$ we choose  $b_i$ with $0\leq b_i<a_i$ such that
$$\sum_{i\in A}b_i+\sum_{i\notin A}a_i=d-1,$$
then the monomial $z=\prod_{i\in A}x_i^{b_i}\prod_{i\notin A}x_i^{a_i}$ satisfies $I:z=P_A$.
\end{Proposition}

\begin{proof}
$\Leftarrow$:  Since $\sum_{i\not\in A} a_i\leq d-1$, it follows
that for any monomial $u\in G(I)$ there exists an integer $j$ with
$j\in A$ such that $x_j$ divides $u$. Therefore $P_A\supset I$.
The condition $\sum_{i=1}^n a_i\geq d-1+|A|$ implies that
$\sum_{i\in A} (a_i-1)+\sum_{i\not\in A} a_i\geq d-1$, which
together with $\sum_{i\not\in A} a_i\leq d-1$ shows that there 
exists  for all $i\in A$ an integer  $b_i$, with $0\leq b_i<a_i$
such that $\sum_{i\in A} b_i+\sum_{i\not\in A} a_i=d-1$. The
monomial  $z=\prod_{i\in A} x_i^{b_i}\prod_{i\not\in A} x_i^{a_i}$
has degree $d-1$. Hence $z\not\in I$, and obviously $I:z\supset
P_A$. If equality holds then $P_A\in\Ass(S/I)$. Assume by
contradiction that $P_A$ is a proper subset of $I:z$. Then there
exists a monomial $u'\in K[x_i:i\not\in A]$ of degree at least
$1$, such that $zu'\in I$. That is,  there exists  a monomial
$u\in G(I)$, $u=\prod_{i=1}^n x_i^{c_i}$, with $u|zu'$. Therefore
$c_i\leq b_i$ for any $i\in A$, since $u'\in K[x_i:i\not\in A]$.
But then
$$d=\deg(u)=\sum_{i=1}^n c_i\leq \sum_{i\in A} b_i+\sum_{i\not\in A}
a_i=\deg(z)=d-1,$$ a contradiction. Therefore $I:z=P_A$, and we are
done.

\implies:  Assume now that $P_A\in\Ass(S/I)$. Then there exists a
monomial $z\in S$, $z\not\in I$, such that $I:z=P_A$. We  first
show that we can choose a monomial $z$ of degree $d-1$ such that
$I:z=P_A$.

Suppose that the monomial $z\not\in I$, with $I:z=P_A$, has degree
$\deg(z)\geq d$ and is of the form $z=\prod_{i=1}^n x_i^{b_i}$. We
observe that there exists an integer $j_0\in [n]$ such that
$b_{j_0}>a_{j_0}$. Indeed, if $b_i\leq a_i$ for all $i\in [n]$,
then $z\in I$, a contradiction.

Since $I:z=P_A$, we have  $x_iz\in I$ for all $i\in A$, and
$x_iz\not\in I$ for all $i\not\in A$. Furthermore, for each $i\in
A$ there exists a monomial $u_i\in G(I)$ such that $u_i$ divides
$x_iz$. This fact together with $z\not\in I$ implies that for all
$i\in A$ the variable $x_i$ appears in $u_i$ with the exponent
$b_i+1$, therefore $b_i<a_i$ for all $i\in A$. In particular, it
follows that $j_0\not\in A$.

We claim that $z/x_{j_0}\not\in I$ and  $I:(z/x_{j_0})=P_A$. The
first assertion follows from that $z\not\in I$ and the fact that
$b_{j_0}-1\geq a_{j_0}$. For the second, we use the general fact
that if $z$ and $z'$ are monomials such that $z'$ divides $z$,
then $I:z'\subset I:z$. Therefore, $I:(z/x_{j_0})\subset P_A$.
Since $b_{j_0}-1\geq a_{j_0}$, then $u_i$ still divides
$x_iz/x_{j_0}$ for all $i\in A$, so $x_i\in I:(z/x_{j_0})$ for all
$i\in A$.  Hence we have the other inclusion $P_A\subset
I:(z/x_{j_0})$, and  our claim is proved.

After a finite number of such reductions, we find  $z\not\in I$ of
degree $d-1$  such that $I:z=P_A$.

In what follows we shall prove that $z$ has the required form of
the statement. As above let $z=\prod_{i=1}^n x_i^{b_i}$ with
$\deg(z)=d-1$, such that $I:z=P_A$. From $I:z=P_A$, it follows
that $zx_i\in I$, for all $i\in A$, and $zx_i\not\in I$, for all
$i\not\in A$. In particular, $b_i+1\leq a_i$ for all $i\in A$, and
$b_i\leq a_i$ for all $i\not\in A$. Since $zx_i\not\in I$ for all
$i\not\in A$,  and since $b_i\leq a_i$, we must  necessarily have
$b_i=a_i$, for all $i\not\in A$. Therefore, we obtain that
$$z=\prod_{i\in A} x_i^{b_i}\prod_{i\not\in A}x_i^{a_i},$$
 with $0\leq b_i<a_i$ for all $i\in A$.
In particular, we have $$\sum_{i\not\in A} a_i\leq\deg(z)=d-1,$$
and $$\sum_{i=1}^n a_i\geq \sum_{i\in A}(b_i+1)+\sum_{i\not\in
A}a_i=d-1+|A|.$$ This concludes the proof of the proposition.
\end{proof}

\begin{Example}
\label{pandurii}
{\em In this example, we compute with the help of the Proposition \ref{radoi}, the set of all associated prime ideals of the ideal of Veronese type $I_{5;3,2,1}\subset K[x_1,x_2,x_3]$. Hence, $d=5$, $a_1=3,a_2=2,a_3=1$ and  $G(I_{5;3,2,1})=\{x_1^3x_2^2, x_1^3x_2x_3, x_1^2x_2^2x_3\}$. Since $a_2+a_3\leq 4$, it follows from the above proposition that $(x_1)=P_{\{1\}}$ is an associated prime ideal, with $I_{5;3,2,1}:x_1x_2^2x_3=(x_1)$, where $z=x_1x_2^2x_3$ and $b_1=1< a_1$. Similarly, one can check that the other associated prime ideals of $I_{5;3,2,1}$ are $(x_2)$, $(x_1,x_2)$, $(x_1,x_3)$, $(x_2,x_3)$. They are obtained as follows: $(x_2)=I_{5;3,2,1}:x_1^3x_3$, $(x_1,x_2)=I_{5;3,2,1}:x_1^2x_2x_3$, $(x_1,x_3)=I_{5;3,2,1}:x_1^2x_2^2$, $(x_2,x_3)=I_{5;3,2,1}:x_1^3x_2$. In particular, $I_{5;3,2,1}$ is equidimensional and is not unmixed. }
\end{Example}

An easy consequence of Proposition \ref{radoi} is the following

\begin{Corollary}
\label{hagi} Let $I\subset S$ be an ideal of Veronese type, and
suppose that  $P_A\in \Ass(S/I)$ for some $A\subset [n]$ with
$|A|=k$. Let $i\in A$,  $j<i$ such that $j\not\in A$  and
$B=(A\setminus \{i\})\union \{j\}$. Then $P_B \in \Ass(S/I)$. In
particular, $(x_1,\ldots,x_k)\in\Ass(S/I)$.
\end{Corollary}
\begin{proof}
Let $A\subset [n]$, such that $|A|=k$ and $P_A\in\Ass(S/I)$. It
follows from Proposition \ref{radoi} that $\sum_{l=1}^n a_l\geq
d-1+k$ and $\sum_{l\not\in A}a_l\leq d-1$. Then, the subset $B$ from the
hypothesis satisfies $|B|=k$ and $$\sum_{l\not\in B} a_l
=\sum_{l\not\in A}a_l-a_j+a_i.$$ Because $j<i$, then $a_i-a_j\leq 0$
and therefore $\sum_{l\not\in B} a_l\leq d-1$. Applying again
Proposition \ref{radoi} we obtain that $P_B\in\Ass(S/I)$. The last
assertion follows after performing at most $k$ times the above
procedure.
\end{proof}

Let us notice that for an ideal of Veronese type
$I=I_{d;a_1,\ldots,a_n}$ with $d=a_1=\ldots =a_i$ for some  $1\leq
i\leq n$, then the maximal graded ideal $\mm$ of $S$ is an
associated prime ideal of $I$ (since $\mm=I:(x_1^{d-1})$). It
follows that if $I$ is unmixed then  $I$ is equidimensional.
According to the Proposition \ref{radoi}, $\mm$ is the only
associated prime of $S/I$ if and only if $a_1=a_2=\ldots=a_n=d$.
Indeed, if for example $a_n<d$, then Proposition \ref{radoi} tells
us that $(x_1,\ldots,x_{n-1})$ is an associated prime ideal of
$S/I$, hence $I$ is not unmixed, a contradiction. Therefore, for
an ideal $I=I_{d;a_1,\ldots,a_n}$ with $d=a_1=\ldots =a_i$ for some  $1\leq i\leq
n$ the following conditions are equivalent: (i) $I$ is unmixed,
(ii) $I$ is Cohen-Macaulay, and (iii)  $I$ is a Veronese ideal.

\begin{Theorem}
\label{vlad} Let $I\subset S$ be a non-zero ideal of Veronese
type. Then the following conditions are equivalent:
\begin{enumerate}
\item[(a)] $I$ is a Cohen-Macaulay ideal,
\item[(b)] $I$ is unmixed and equidimensional,
\item[(c)] $I$ is unmixed.
\end{enumerate}
\end{Theorem}
\begin{proof}
Following the discussion above theorem it is enough to prove for an
ideal of Veronese type, $I$, with $d>a_1\geq\ldots\geq a_n\geq 1$.
The implications $(a)\Longrightarrow (b)$, $(b)\Longrightarrow (c)$
are obvious.

To complete the proof of the theorem it remains to prove
$(c)\Longrightarrow (a)$. Herzog and Hibi proved (see \cite[Theorem
3.2.]{HH1}) that a polymatroidal ideal $I$ is a Cohen-Macaulay ideal
if and only if  $I$ is either a principal ideal, a Veronese ideal,
or a squarefree Veronese ideal. From the discussions preceding the
theorem it follows that we only have to show that $I$ is unmixed
implies that $I$ is either principal or squarefree Veronese.

First we shall show that $I$ is equidimensional. Assume by contrary
that $I$ is not equidimensional. Then $\sqrt{I}$ is not
equidimensional and therefore there exist $A,B\subset [n]$ of
cardinality $k$, respectively $l$ such that $k<l$ and $P_A,P_B$ are
minimal primes of $\sqrt{I}$, and consequently of $I$. It follows
from Corollary \ref{hagi} that $(x_1,\ldots,x_k)\subset
(x_1,\ldots,x_l)$ are both associated prime ideals of $I$, a
contradiction since $I$ is unmixed. Therefore $I$ is
equidimensional. Hence, all the associated prime ideals of $S/I$ are
minimal over $I$ and have the same height, which we denote it by
$k$, with $1\leq k\leq n$. Now we have two cases.
\begin{enumerate}
\item[Case 1:] $a_1+\ldots+a_n=d$.
\end{enumerate}
This implies that $I$ is principal, generated by $x_1^{a_1}\cdots
x_n^{a_n}$, hence Cohen-Macaulay, so we are done.
\begin{enumerate}
\item[Case 2:] $a_1+\ldots+a_n>d$.
\end{enumerate}
Since all the associated primes of $S/I$ have the same height $k$,
it follows from Corollary \ref{hagi} that, in particular,
$(x_1,\ldots,x_k)$ is an associated prime of $S/I$. First we notice
that $k<n$. Indeed, if $k=n$, then by Proposition \ref{radoi}, it
would follow that $a_1+\ldots+a_n\geq d-1+n$. In particular,
$a_1+\ldots+a_n>d-1+(n-1)$ and together with $a_n\leq d-1$, from the
hypotheses, we obtain, via the same proposition, that $(x_1,\ldots,
x_{n-1})$ is an associated prime ideal, a contradiction, since $I$
is unmixed.

Hence $k<n$ and we have that $a_{k+1}+\ldots+a_{n}\leq d-1$ and
$a_1+\ldots+a_n\geq d-1+k$. Because $I$ is unmixed, then for any
subset $A\subset [n]$ with $|A|\neq k$, $P_A\not\in\Ass(S/I)$.
Therefore $(x_1,\ldots,x_{k+1})\not\in\Ass(S/I)$ and since
$a_{k+2}+\ldots+a_n<a_{k+1}+\ldots+a_n\leq d-1$, it follows from
Proposition \ref{radoi} that $a_1+\ldots+a_n<d-1+k+1=d+k$. This
inequality together with $a_1+\ldots+a_n\geq d-1+k$ imply that
$a_1+\ldots+a_n=d-1+k$. Since we are in the case $a_1+\ldots+a_n>d$,
we necessarily have $k>1$.

We conclude this case by showing that $a_i=1$ for all $i$ with
$1\leq i\leq n$, which means that $I$ is squarefree Veronese ideal,
and hence a Cohen-Macaulay ideal. Suppose, by contradiction, that
there exists $i\in [n]$ such that $a_i\geq 2$, hence in particular
$a_1\geq 2$. Then, because $k>1$ and $I$ is unmixed, we have
$(x_1,\ldots,x_{k-1})\not\in\Ass(S/I)$. Since
$a_1+\ldots+a_n=d+k-1>d-1+k-1$, applying again Proposition
\ref{radoi} we obtain that $a_k+\ldots+a_n>d-1$. Suming up this
inequality with the inequality $a_1+\ldots+a_{k-1}\geq k$ (true,
because we supposed that $a_1\geq 2$) we obtain that
$a_1+\ldots+a_n> d-k+1$, a contradiction. Therefore, $a_i=1$ for all
$i\in [n]$, and we are done.

In conclusion, from the two cases analyzed we obtain that if $I$ is
unmixed then $I$ is Cohen-Macaulay, so the theorem is proved.
\end{proof}

\end{document}